\documentclass[a4paper,10pt]{amsart}

\usepackage{amstext}
\usepackage{amsmath}
\usepackage{ifthen}
\usepackage{amscd}
\usepackage{amssymb}
\usepackage[english]{babel}
\usepackage[T1]{fontenc}
\usepackage[utf8]{inputenc}
\usepackage[all]{xy}
\usepackage{graphicx}
\usepackage{enumerate}

\newenvironment{demo}[1][]{\ifthenelse{\equal{#1}{}}{\noindent \textbf{Proof :\\ \indent}}{\noindent \textbf{Proof #1 :\\ \indent}}}{$\square$\\}

\newtheorem{thm}{Theorem}[section]
\newtheorem{conj}{Conjecture}[section]
\newtheorem{cor}{Corollary}[section]
\newtheorem{lem}{Lemma}[section]
\newtheorem{prop}{Proposition}[section]

\newtheorem{rem}{Remark}[section]

\newcommand{\CC}{\mathbb C}

\newcommand{\ZZ}{\mathbb Z}

\newcommand{\DD}{\mathbb D}
\newcommand{\PP}{\mathbb P}

\newcommand{\kod}[1]{\kappa(#1)}
\newcommand{\gf}[1]{\pi_1(#1)}

\newcommand{\To}{\longrightarrow}
\newcommand{\abs}[1]{\left\vert#1\right\vert}

\newcommand{\set}[1]{\left\{#1\right\}}

\newcommand{\dimm}[1]{\mathrm{dim}(#1)}

\newcommand{\wtx}{\widetilde{X}}
\newcommand{\formL}[1]{H^0_{(2)}(\wtx,\Omega^{#1}_{\wtx})}

\title{Non algebraicity of universal covers of K\"ahler surfaces}
\author{Beno\^it \textsc{Claudon}}
\address{Institut Fourier - UMR 5582 - 100, rue des Maths B.P. 74 38402 Saint-Martin d'H\`eres, France}
\email{Benoit.Claudon@ujf-grenoble.fr}

\begin{document}

\maketitle

\begin{abstract}
In this short note, we prove the following fact : if the universal cover of a K\"ahler compact surface is biholomorphic to an affine variety, then, up to some finite \'etale cover, the surface is a torus.
\end{abstract}

\section{Introduction}

In his book \cite{Sh74}, I. Shafarevich splits the classification of universal covers of projective manifolds into two types. The type I corresponds to manifolds having a finite fundamental group and their universal covers still belong to the class of projective manifolds. The type II is (according to him) somehow more delicate to define but in this case the fundamental group is infinite and the universal cover is very far from being a projective manifold\footnote{the exact quotation is the following : \emph{Type II is hard to define at present, beyond the fact the fundamental group is infinite. In this case, the universal cover is a ``very big'' complex manifold, very far from being a projective or complete variety.}}.

The following question seems then to be natural : how far from being projective this universal cover can be ? Since there is the obvious example of complex tori, we can wonder if the universal cover of a projective manifold (or more generally of a compact K\"ahler manifold) can appear as a quasi-projective variety (or a quasi-K\"ahler manifold in the K\"ahler setting). The case of curves is easily settled : the elliptic curves are the only one having in the same time an infinite fundamental group and a universal cover which can be compactified.

In this article, we study the surface case and show the following result. 

\begin{thm}\label{main}
Let $X$ be a compact K\"ahler surface whith infinite fundamental group and whose universal cover is quasi-K\"ahler (\emph{i.e.} can be realized as a Zariski open subset of a compact K\"ahler surface). Then, $X$ is one of the following :
\begin{enumerate}[(1)]
\item $X$ is geometrically ruled on an elliptic curve and $\wtx\simeq\CC\times\PP^1$,
\item $X$ admits then a finite \'etale cover which is a torus and the universal cover is $\CC^2$.
\end{enumerate}
\end{thm}
We can specialize it to the affine case.
\begin{cor}
Let $X$ be a compact K\"ahler surface whose universal cover $\wtx$ is (biholomorphic to) an affine variety. Then, $\wtx\simeq\CC^2$ and, up to some finite \'etale cover, $X$ is a two-dimensional torus.
\end{cor}

\noindent It is quite natural to formulate the following hope :
\begin{conj}\label{conjecture affine}
Let $X$ be a compact K\"ahler manifold whose universal cover is (biholomorphic to) an affine variety. Up to some finite \'etale cover, $X$ is a torus.
\end{conj}
This conjecture should be related to the Iitaka conjecture (and can be seen as a strengthening of the latter).
\begin{conj}\label{Iitaka}
A compact K\"ahler manifold covered by $\CC^n$ has a finite \'etale cover which is a torus.
\end{conj}
Note that the conjecture \ref{Iitaka} holds up to dimension 3 ; see \cite{CZ04}.\\

To conclude this introduction, let us point out the fact that the K\"ahler hypothesis is essential to ensure the validity of theorem \ref{main}. There exist several examples of non K\"ahler compact surfaces whose universal covers appear as quasi-projective varieties. Among them, we can quote the following examples :
\begin{enumerate}
\item Hopf surfaces: the universal cover is then $\CC^2\backslash(0)$.
\item Kodaira surfaces: these surfaces are elliptic submersion on an elliptic curve. Their universal covers are isomorphic to $\CC^2$ and their fundamental group are solvable groups which sit in an exact sequence
$$1\To\ZZ^2\To\gf{X}\To\ZZ^2\To1.$$
Les us note that these groups are not almost abelian (not even almost nilpotent).
\end{enumerate}
A detailed description of these examples can be found in \cite{BPV}.

\section{The non-general type case}

In this section, we study the case where the surface $X$ satisfies $\kod{X}\le1$. For sake of brievity, let us denote (in the rest of this article) by $(H)$ the following hypothesis :
\begin{center}
$(H)$ : $\gf{X}$ is infinite and $\wtx$ is quasi-K\"ahler.  
\end{center}

\noindent We begin this section by the following observation : the hypothesis $(H)$ implies minimality.

\begin{lem}\label{qK implique minimal}
If a surface $X$ satisfies the hypothesis $(H)$, it is minimal : $X$ does not contain a $(-1)$-curve.
\end{lem}
\begin{demo}
Since $\wtx$ is assumed to be quasi-Kähler, there exists a compact surface $\hat{X}$ such that $\wtx\hookrightarrow\hat{X}$ is an open subset. If $C$ is a $(-1)$-curve contained in $X$, $\pi^{-1}(C)$ can then be written as a disjoint union of rational curves in $\wtx$ and hence in $\hat{X}$ :
$$\pi^{-1}(C)=\bigcup_{i\in\pi_1(S)}C_i.$$
Since being a $(-1)$-curve is a local property of the embedding of a curve, each curve $C_i$ has self-intersection -1 in $\hat{X}$. Here we have a contradiction : a compact surface cannot contain infinitely many \emph{disjoint} $(-1)$-curves.
\end{demo}
\begin{rem}\label{exemple Zariski}
A compact surface cannot contain infinitely many \emph{disjoint} $(-1)$-curves because the second Betti number strictly decreases at each blow-down. It should be noted that a rational surface can contain infinitely many $(-1)$-curve. For instance, the blow-up of $\PP^2$ along the base locus of a pencil of cubics has this property.
\end{rem}

The main point is now to understand which elliptic surfaces can in the same time fullfill the hypothesis \emph{(H)}.
\begin{prop}\label{surface elliptique-qK}
Let $X\To C$ be an elliptic surface over a curve $C$. If $X$ satisfies the hypothesis \emph{(H)}, then the universal cover of $X$ is either $\CC^2$ (and $X$ is covered by a torus) or $\PP^1\times\CC$. In particular, in both cases we have $\kod{X}\leq 0$.
\end{prop}
\begin{demo}
We can assume (replacing $S$ by a finite \'etale cover) that the fibration induces an exact sequence (see for instance \cite[appendix C]{Ca98}):
\begin{equation}\label{suite pi1}
1\To\pi_1(E)_X=\mathrm{Im}\left(\pi_1(E)\To\pi_1(X)\right)\To \pi_1(X)\To\pi_1(C)\To 1.
\end{equation}
The universal cover of $X$ being quasi-Kähler do not admit any non constant bounded holomorphic function and this implies that $g(C)\le1$. If the fibration has only multiple smooth fibres (of type $mI_0$), some finite \'etale cover of $X$ will become a smooth elliptic fibration. We can then apply the proposition \ref{sub elliptique} below to conclude that $\wtx\simeq\CC\times\tilde{C}$.

Let us assume now that the fibration has a genuine singular fiber (over $y\in C$). We can then infer that the quotient $\pi_1(E)_X$ is finite (as in \cite{GS85}) and that $g(C)=1$ ($\wtx$ is assumed to be non compact). This eaxctly means that the induced fibration $\wtx\To \CC$ is a proper map. The family of fibres of this fibration consists in an analytic family of compact cycles of $\wtx$ and then it induces a family of cycles of $\hat{X}$. The cycles space $\mathcal{C}(\hat{X})$ (constructed in \cite{Ba75}) has only compact irreducible components and the fibration can then be extended to a (proper) fibration $\hat{X}\To\hat{C}=\PP^1$. But this fibration has one singular fibre for each point $\tilde{y}\in\hat{C}$ lying over $y\in C$ ; since it consists of infinitely many points, we get a fibration of a compact surface with infinitely many singular fibers. This is clearly impossible and the proposition is proved.
\end{demo}

\noindent In the course of the preceeding proof, we used the following fact.
\begin{lem}\label{sub elliptique}
Let $f:X\To Y$ a submersion between \emph{compact} complex manifolds whose fibers are genus 1 curves. If moreover $X$ is supposed to be K\"ahler, then $f$ is isotrivial and $\wtx$ is biholomorphic to $\CC\times\widetilde{Y}$.
\end{lem}

\begin{demo}
The $j$ function
$$j:\left\{\begin{array}{rcl}Y & \To & \CC\\
y&\mapsto & j\left(f^{-1}(y)\right)
\end{array}\right.
$$
is well-defined: use local (on $Y$) analytic sections of $f$ to have genuine elliptic curves and compute $j$. Since $f$ is a submersion, $j$ is holomorphic on $Y$ hence constant. This exactly means that all fibers are isomorphic to a fixed elliptic curve $E$; in particular, up to a finite \'etale cover, $f:X\To Y$ is a principal bundle of fiber $E$. From the fact that $X$ is K\"ahler, we infer that $q(X)=q(Y)+1$ or, equivalently, that the morphism
$$i_*:\gf{E}\To\gf{X}$$
induced by $i:E\hookrightarrow X$ is injective. We can then consider the Albanese tori of $X$ and $Y$; the map
$$\alpha(f):\mathrm{Alb}(X)\To\mathrm{Alb}(Y)$$
is an elliptic bundle and induced clearly an isomorphism
$$\widetilde{\mathrm{Alb}(X)}\simeq\widetilde{\mathrm{Alb}(Y)}\times\CC.$$
Pulling back this splitting to $X$ yields the proposition.
\end{demo}

\begin{rem}\label{hopf}
The K\"ahler hypothesis is certainly not superfluous. If $S\To\PP^1$ is an elliptic Hopf surface, the morphism
$$\ZZ^2=\gf{E}\stackrel{i_*}{\To}\gf{S}=\ZZ$$
is absolutely not injective and the universal cover of $S$ is $\CC^2\backslash\set{0}$, \emph{i.e.} a $\CC^*$-bundle over $\PP^1$.
\end{rem}

We can now prove the main theorem under the assumption that the surface $X$ is not of general type.\\

\begin{demo}[of the theorem \ref{main}, case $\kod{X}<2$]
Being minimal (lemma \ref{qK implique minimal}) and not of general type, the surface $X$ belongs to one of the following class :
\begin{description}
\item[\mathversion{bold}$\kod{X}=-\infty$\mathversion{normal}] $X$ is geometrically ruled on a curve $C$ ($X=\PP(E)$ for some rank 2 vector bundle $E$ on $C$),
\item[\mathversion{bold}$\kod{X}=0$\mathversion{normal}] up to a finite \'etale cover, $X$ is a torus or a $K3$ surface,
\item[\mathversion{bold}$\kod{S}=1$\mathversion{normal}] $X$ is elliptic over a curve $C$.
\end{description}
In the first case, we can easily see that $g(C)=1$ and, the bundle $E$ becoming trivial on $\tilde{C}=\CC$, that $\wtx\simeq \PP^1\times\CC$ (which is actually a quasi-projective variety). In the $\kod{X}=0$ case, the fact that $\pi_1(X)$ is infinite implies that $X$ is a torus. The last case is ruled out by proposition \ref{surface elliptique-qK}.
\end{demo}

\section{$L^2$ canonical forms and the general type case}

In this final section, we have to prove that a (projective) surface cannot in the same time be of general type and satisfy the hypothesis \emph{H}. To achieve this goal, we resort to $L^2$ theory on $\wtx$.

\subsection{Review of $L^2$ index theory}

In this paragraph, we state a particular case of the $L^2$ index theorem of Atiyah \cite{At76}; for more on this subject, see \cite{At76} or the books \cite{Lu02} and \cite{Roe98}.

On $\wtx$, consider the spaces of $(p,q)$ harmonic forms which are $L^2$ with respect to $\tilde{\omega}$ :
$$\mathcal{H}^{p,q}(\wtx)=\set{\alpha\, (p,q)-\textrm{form}\vert \,\Delta\alpha=0,\,\int_{\wtx}\abs{\alpha}_{\tilde{\omega}}dV_{\tilde{\omega}<\infty}}.$$
The case $(p,0)$ has the following interpretation; the K\"ahler identities show that $(p,0)$ harmonic forms are indeed holomorphic and theses spaces can be described as follows:
$$\mathcal{H}^{p,0}(\wtx)=H^0_{(2)}(\Omega^p_{\wtx})=
\set{\alpha\in H^0(\Omega^p_{\wtx})\vert\,\int_{\wtx}\alpha\wedge\overline{\alpha}\wedge\tilde{\omega}^{n-p}<\infty}.$$
In particular, for $p=n$, we have the Bergman space and the integrability condition is independant of the metric $\tilde{\omega}$.\\

These spaces of $L^2$ harmonic forms are Hilbert spaces which are moreover endowed with an isometric action of $\Gamma=\gf{X}$. Using this action (and the compactness of $X$), it is possible to define a renormalized dimension function (which takes \emph{a priori} non-negative real values); we will use the following notations (we omit the additional subscript $\Gamma$):
$$h^{p,q}_{(2)}(\wtx):=\mathrm{dim}_{\Gamma}\left(\mathcal{H}^{p,q}(\wtx)\right).$$
This notion of dimension enjoy some usual properties; we will need one of them in the sequel.
\begin{prop}\label{G-dimension}
If $h^{p,q}_{(2)}(\wtx)>0$ and if $\Gamma$ is infinite, the space $\mathcal{H}^{p,q}(\wtx)$ is infinite dimensional (in the usual sense).
\end{prop}
\begin{demo}
If $\mathcal{H}^{p,q}(\wtx)$ were finite dimensional (but non zero), we could consider the function:
$$k(x)=\sum_{j=1}^N \abs{\alpha_j(x)}_{\tilde{\omega}}^2\quad(x\in\wtx)$$
where $(\alpha_j)_{j=1..N}$ is an orthonormal basis of $\mathcal{H}^{p,q}(\wtx)$. Since this expression is independant of the chosen basis, it is $\Gamma$-invariant and the computation of its integral on $\wtx$ gives:
$$\int_{\wtx}k(x)dV_{\tilde{\omega}}=\sum_{j=1..N} \int_{\wtx}\abs{\alpha_j}_{\tilde{\omega}}^2dV_{\tilde{\omega}}=N.$$
Since we assume $\Gamma$ to be infinite, the first integral is in fact infinite ($k$ is $\Gamma$-invariant and not identically zero) and we get a contradiction.
\end{demo}

In this setting, the $L^2$ index theorem can be stated in the following form.
\begin{thm}[M. F. Atiyah, \cite{At76}]\label{indice L2}
The $L^2$ Euler characteristic of $\wtx$ is the same as the Euler characteristic of $X$:
$$\chi_{(2)}(\mathcal{O}_{\wtx}):=\sum_{p=0}^n(-1)^p h^{p,0}_{(2)}(\wtx)=\chi(\mathcal{O}_X).$$
\end{thm}

\subsection{Existence of $L^2$ canonical forms and conclusion}

We apply the $L^2$ index theorem to get $L^2$ canonical forms on universal covers of certain projective manifolds of general type.

\begin{lem}\label{existence L2}
Let $X$ be a projective manifold of general type. Then there exists an $L^2$ canonical form on $\widetilde{X}$ in the following situations:
\begin{enumerate}[(i)]
\item $\dimm{X}=2$ and $\gf{X}$ is infinite.
\item $\dimm{X}=3$, $X$ has a smooth minimal model and $\wtx$ is not covered by compact surfaces.
\end{enumerate}
\end{lem}
\begin{demo}
Let us consider first the surface case. Since everything is invariant under birational morphism, we can assume that the surface $X$ is minimal: $K_X$ is then \emph{big} and \emph{nef}. Applying $L^2$-index theorem and Miyaoka-Yau inequality, this yields:
$$h^0_{(2)}(K_{\widetilde{X}})-h^0_{(2)}(\Omega^1_{\widetilde{X}})=\chi_{(2)}(\mathcal{O}_{\widetilde{X}})=\chi(\mathcal{O}_X)=\frac{c_1^2+c_2}{12}\geq\frac{1}{9}c_1^2>0.$$
We use here the vanishing of $h^0_{(2)}(\mathcal{O}_{\wtx})$: being $L^2$, an holomorphic function tends to zero at infinity; it has to be constant (maximum principle) and then zero since $\wtx$ has infinite volume. The positivity of $c_1^2$ ($K_X$ being nef and big) implies then the positivity of $h^0_{(2)}(K_{\widetilde{X}})$.

The three-dimensional case is based on the same arguments. Since $\wtx$ is not covered by compact surfaces, Gromov Theorem \ref{Gromov} below implies the vanishing of $h^0_{(2)}(\Omega^1)$ and the same computation as above gives:
$$h^0_{(2)}(\Omega^2_{\widetilde{X}})-h^0_{(2)}(K_{\widetilde{X}})=\frac{c_1c_2}{24}\leq\frac{1}{64}c_1^3<0.$$
Here we use the fact that we can assume $X$ smooth and minimal so as to apply the Miyaoka-Yau inequality. This proves the second statement.
\end{demo}

Let us recall the statement of the theorem quoted in the previous proof.
\begin{thm}[M. Gromov, \cite{G91}]\label{Gromov}
Let $X$ be a compact K\"ahler manifold with non vanishing first $L^2$ Betti number $h^0_{(2)}(\Omega^1_{\wtx})$ and with infinite fundamental group. Then there exists a \emph{proper} (surjective) map
$$f:\wtx\To \DD$$
such that $H^0_{(2)}(\Omega^1_{\wtx})$ is induced by $f$; in particular, $\wtx$ is covered by compact hypersurfaces.
\end{thm}

\begin{rem}\label{green-laz}
The smoothness of the minimal model cannot be ommited in lemma \ref{existence L2}. To see this consider the example described in \cite{EL97}. Let $\pi:C\To E$ be a double cover of an elliptic curve with $g(C)\ge 2$ and denote by $\tau$ the involution of $\pi$. If $A=E\times E\times E$ and $X$ is a desingularization of the quotient
$$Y=C\times C\times C/<\tau\times\tau\times\tau>,$$
$X$ has a generically finite map $X\To A$ (its Albanese variety). It is obvious that $X$ is of general type; nevertheless its universal cover has no non-zero $L^2$ canonical form. Actually, it is proven in \cite{EL97} that $\chi(\mathcal{O}_X)=0$; being generically finite on an abelian variety, $X$ satisfies the following $L^2$ vanishings conditions :
$$\formL{j}=0,\,j=0,1,2.$$
The $L^2$-index theorem is reduced to the equality we looked for
$$h^0_{(2)}(K_{\wtx})=-\chi_{(2)}(\mathcal{O}_{\wtx})=-\chi(\mathcal{O}_X)=0.$$
It should be pointed out that the minimal model of $X$ is nothing but $Y$ which is (fortunately) not smooth.
\end{rem}

We can now finish the\\

\begin{demo}[of the theorem \ref{main}]
Indeed, if $X$ is surface of general type, its universal cover carries some $L^2$ canonical forms (lemma \ref{existence L2}). On the other hand, if $X$ satisfies the hypothesis $(H)$, we can consider some compactification $\wtx\hookrightarrow\hat{X}$ of $\wtx$. It is now a well known fact the $L^2$ canonical forms extend through the divisor $D=\hat{X}\backslash\wtx$. Indeed, choosing coordinates adapted to $D$ near a smooth point of $D$ and using Taylor expansion to compute the $L^2$ norm, we see easily that the $L^2$ condition prevent the canonical forms to have poles on $D$ (here we use the important fact that the $L^2$ condition does not involve the metric $\tilde{\omega}$; this is true only for canonical forms).\\
We get thus an injection:
$$H^0_{(2)}(\wtx,K_{\wtx})\hookrightarrow H^0(\hat{X},K_{\hat X}).$$
At this point, we have a contradiction: $H^0(\hat{X},K_{\hat X})$ is finite dimensional (compactness of $\hat{X}$) whereas $H^0_{(2)}(\wtx,K_{\wtx})$ is infinite dimensional (proposition \ref{G-dimension}). This concludes the proof of the main theorem.
\end{demo}

\bibliographystyle{amsalpha}
\bibliography{myref}

\end{document}